\documentclass[reqno]{amsart}
\usepackage{amsmath,amssymb,amsfonts,hyperref}

\newcommand{\FR}{{\rm FR}}
\newcommand{\ST}{{\rm ST}}

\newcommand{\TTH}{{\rm TH}}

\newcommand{\STAB}{{\rm STAB}}

\newcommand{\ol}{\overline}

\newcommand{\PP}{{\mathcal P}}
\newcommand{\AAA}{{\mathcal A}}
\newcommand{\tL}{\widetilde L}
\newcommand{\MM}{{\mathcal M}}
\newcommand{\MMV}{{\mathcal M}_{+,V}}
\newcommand{\SSS}{{\mathcal S}}

\newcommand{\oR}{\mathbb R}

\newcommand{\bI}{{\bf I}}

\newcommand{\bO}{\mathbf 0}
\newcommand{\xx}{\left(\hspace*{-1.7mm}\begin{array}{c} 1 \\ x \end{array}\hspace*{-1.7mm} \right)}
\newcommand{\conv}{{\rm conv}}   

\newtheorem{theorem}{\bf Theorem}[section]

\newtheorem{lemma}[theorem]{\bf Lemma}
\newtheorem{example}[theorem]{\bf Example}

\newtheorem{proposition}[theorem]{\bf Proposition}

\theoremstyle{remark}

\begin{document}

\title{Block-diagonal semidefinite programming hierarchies for 0/1
programming}

\author{Neboj\v{s}a Gvozdenovi\'c} 
\address{N.~Gvozdenovi\'c, The
  Faculty of Economics, University of Novi Sad, Segedinski put 9-11,
  24000 Subotica, Serbia} \email{nebojsa.gvozdenovic@gmail.com}

\author{Monique Laurent}
\address{M.~Laurent, Centrum voor Wiskunde en Informatica (CWI),
Kruislaan 413, 1098 SJ Amsterdam, The Netherlands}
\email{monique@cwi.nl}

\author{Frank Vallentin}
\address{F.~Vallentin, Centrum voor Wiskunde en Informatica (CWI),
Kruislaan 413, 1098 SJ Amsterdam, The Netherlands}
\email{f.vallentin@cwi.nl}

\thanks{This work was supported by the Netherlands Organization for
Scientific Research under grant NWO 639.032.203. The third author was
partially supported by the Deutsche Forschungsgemeinschaft (DFG) under
grant SCHU 1503/4.}



\date{August 19, 2008}

\begin{abstract}
  Lov\'asz and Schrijver, and later Lasserre, proposed hierarchies of
  semidefinite programming relaxations for general 0/1 linear
  programming problems.  In this paper these two constructions are
  revisited and two new, block-diagonal hierarchies are proposed. They
  have the advantage of being computationally less costly while
  being at least as strong as the Lov\'asz-Schrijver hierarchy. Our construction is
  applied to the stable set problem and experimental results for Paley
  graphs are reported.
\end{abstract}

\maketitle

\markboth{N.~Gvozdenovi\'c, M.~Laurent, F.~Vallentin}{Block-diagonal
semidefinite programming hierarchies}

\vspace*{-1cm}

\section{Introduction}

A basic approach in combinatorial optimization consists of formulating
the problem at hand as a 0/1 linear programming problem, typically of
the form
\begin{equation*}
\max\left\{ c^T x \mid Ax \leq b,\; x \in \{0,1\}^n\right\},
\end{equation*}
where $c\in\oR^n$, $b\in \oR^m$ and $A\in \oR^{m\times n}$.  Then the
task is to find an efficiently computable outer approximation of the
polytope $P$, defined as the convex hull of the 0/1 solutions to $Ax
\leq b$.

On the one hand, extensive research has been done for finding
(partial) linear inequality descriptions for many polyhedra arising
from specific combinatorial optimization problems. On the other hand,
researchers focused on developing general purpose methods for
arbitrary 0/1 linear programming problems. Here let us mention the
method of Gomory for generating cuts strengthening the initial linear
relaxation $\{x \in \oR^n \mid Ax \leq b\}$ of $P$ and its various
extensions for generating strong cutting planes (see e.g.\
\cite{NW,W}),  the lift-and-project method \cite{BCC93},
the reformulation-linearization technique \cite{SA90},
the matrix-cut method of Lov\'asz and Schrijver \cite{LS91}, and the
sums of squares and moment method of Lasserre \cite{Las01b}. Some of
these methods are compared in \cite{Lau03}; see also \cite{LR05}. A
common feature of the methods of Lov\'asz-Schrijver and of Lasserre is
that they consider hierarchies involving semidefinite relaxations of
$P$: Convex sets $Q_t$ ($t = 1, \ldots, n+1$) are constructed which
can be described by semidefinite conditions and which form a hierarchy
of increasingly stronger relaxations:
\begin{equation*}
\{x \in \oR^n \mid Ax \leq b\} \supseteq Q_1 \supseteq Q_2 \supseteq \ldots
\supseteq Q_{n+1} = P.
\end{equation*}
The two hierarchies are related; it is shown in \cite{Lau03} that the
hierarchy of Lasserre refines the hierarchy of Lov\'asz-Schrijver.

In this note we revisit these hierarchies and propose two new
ones, which differ in the way of encoding the linear constraints defining the starting linear relaxation of $P$.
Moreover one of them (introduced in Section \ref{secvariation})
can also be defined when the starting relaxation of $P$ is an arbitrary convex body, 
as is the case for the 
Lov\'asz-Schrijver construction.  The new hierarchies are nested between the Lasserre and Lov\'asz-Schrijver hierarchies, but they are less costly to compute. So
they are especially well-suited for implementations. E.g., at given order $t$,
the new hierarchy from 
Section~\ref{secvariation} 
involves $1/(t+1)! n^{t+1} +O(n^{t})$ variables compared to
$2^{t-2}n^{t+1} +O(n^{t})$ variables for the Lov\'asz-Schrijver
hierarchy and to $O(n^{2t})$ variables for the Lasserre hierarchy.
The new hierarchies can be seen as a variation of the Lasserre
hierarchy, where one replaces a large matrix of order $O(n^t)$ by
smaller blocks of order $n+1$ arising by block-diagonalizing suitably
defined principal submatrices of the original large matrix. The
motivation for considering block matrices is that it is
computationally easier to solve a semidefinite program involving many
small blocks rather than one large matrix. Most currently available
interior-point algorithms for semidefinite programming are indeed
designed to exploit block-diagonal matrices.  While the hierarchy of
Lov\'asz and Schrijver is originally defined recursively, we give an
explicit description obtained by ``unfolding'' the recursion.  In this
way, the connection to the new hierarchies becomes transparent (see
Section \ref{sechierarchy} for details).

When applied to the stable set problem, our new construction gives a block-diagonal hierarchy whose first two steps were already used in the literature. 
The first order relaxation
gives the Lov\'asz theta number and the second order one
 gives parameters considered in \cite{KP07,Lau07}
for the stable set problem and in \cite{GL1,GL2} for the coloring
problem. In these applications the computational advantage of the new
hierarchy was of crucial importance.

\subsection*{Contents of the paper} 

In Section \ref{sechierarchy} we first briefly introduce the
constructions of Lov\'asz-Schrijver and of Lasserre. Then we give the
new construction and show how to derive more compact formulations by
block-diagonalization. In Section \ref{secstable} we apply it to the
stable set problem and in Section \ref{secPaley} we 
present some computational results illustrating
the behavior of the new hierarchy for approximating the stability
number of Paley graphs.

\subsection*{Notation} 

Given a finite set $V$, we denote the collection of all subsets of $V$
by $\PP(V)$.  Given a non-negative integer $r$, set $\PP_r(V):=\{I\in
\PP(V)\mid |I|\le r\}$ and $\PP_{=r}(V):=\{I\in \PP(V)\mid |I|= r\}$. By $\bO$ we denote the empty set. Sometimes we identify
$\PP_{=1}(V)$ with $V$, i.e., we write $i$ instead of
$\{i\}$. Furthermore, we sometimes write $ij$ instead of $\{i,j\}$ and
$ijk$ instead of $\{i,j,k\}$, etc. The standard unit vectors in
$\oR^{\PP_1(V)}$ are denoted by $e_{\bO}$, $e_i$ for $i\in V$.

\section{Semidefinite programming hierarchies}
\label{sechierarchy}

Suppose we are given a convex cone $K$ contained in the homogenized
unit cube 
$\{x\in\oR^{\PP_1(V)} \mid 0 \le x_i\le x_\bO \ (i\in V)\}$.
Set
\begin{equation*}
\begin{split}
& P_K:=\conv\left\{x\in\{0,1\}^V\mid \xx \in K\right\},\\
& C_K:=\oR_+\left\{\xx\in K \text{ with } x\in\{0,1\}^V\right\}.
\end{split}
\end{equation*}
The general objective is to find the linear inequality description of
the polytope $P_K$ or, equivalently, of the cone $C_K$.  In
Section~\ref{secLS} we recall the construction of Lov\'asz-Schrijver
which applies to any convex cone $K$. While the original construction
is recursive we propose an explicit semidefinite programming
reformulation. In Section~\ref{secL} we recall the construction of
Lasserre which applies to the case when $K$ is represented by
polynomial inequalities. Here we focus on polyhedral cones $K$ of the form
\begin{equation}\label{coneK}
K= \{x\in \oR^{\PP_1(V)}\mid a_l^Tx\ge 0 \ (l=1,\ldots,m)\},
\end{equation}
where $a_1,\ldots,a_m\in \oR^{\PP_1(V)}$.  In Section~\ref{secnew} we
introduce our new construction, which can be seen as a variation of the
previous methods. We discuss two new hierarchies.  The first one
applies to polyhedral cones $K$ as in \eqref{coneK} and is more
economical than the Lasserre hierarchy while still refining the
Lov\'asz-Schrijver hierarchy.  The second one applies to any convex
cone $K$ and can be seen as a non-recursive analogue of the
Lov\'asz-Schrijver hierarchy having a more compact and explicit
formulation.

\subsection{The Lov\'asz-Schrijver hierarchy}
\label{secLS}

In this section we recall basic facts about the Lov\'asz-Schrijver
hierarchy. For proofs and more details we refer to~\cite{LS91}.
Set
\begin{equation*}
\MMV:=\{Y\in \oR^{\PP_1(V)\times \PP_1(V)} \mid Y\succeq 0,\
Y_{ii}=Y_{\bO i}\ (i\in V)\},
\end{equation*}
where ``$\succeq 0$'' stands for ``is positive semidefinite''. For a
convex cone $K\subseteq \oR^{\PP_1(V)}$ define
\begin{equation*}
\begin{split}
& \MM_+(K):= \{Y \in \MMV \mid Ye_i \in K,\ Y(e_\bO-e_i)\in K\ (i\in V)\}, \\
& N_+(K):=\{Ye_\bO \mid Y\in \MM_+(K)\}.
\end{split}
\end{equation*}
The {\em $t$-th iterate of the Lov\'asz-Schrijver hierarchy} is
$N_+^t(K):=N_+(N_+^{t-1}(K))$ for $t\ge 1$, where
$N_+^1(K):=N_+(K)$ and $N^0_+(K):=K$. It lies between $K$ and $C_K$ and
$N^{t+1}_+(K) \subseteq N^t_+(K)$.
We have $N_+^t(K)=C_K$ for $t = |V|$.  Moreover, for any fixed $t$,
if one can optimize over $K$ in polynomial time then the same holds
for $N_+^t(K)$. 

 In the following proposition we ``unfold'' the
recursive definition of $N_+^t(K)$ and give an explicit semidefinite
programming formulation. Its proof is straightforward and thus omitted.

\begin{proposition}
\label{propeqNt}
A vector $x \in \oR^{\PP_1(V)}$ lies in $N_+^t(K)$ if and only if there
exist a matrix $Y \in \MMV$ and matrices $Y^{\sigma_1, \ldots,
\sigma_{s}}_{i_1, \ldots, i_s} \in \MMV$, with $s = 1, \ldots, t-1$,
$i_1, \ldots, i_s \in V$ and $\sigma_1, \ldots, \sigma_s \in \{\pm
1\}$ satisfying the following conditions:
\begin{enumerate}
\item[(a)] $x = Ye_\bO$.
\item[(b)] For all $s = 1, \ldots, t-1$, $i_1, \ldots, i_s \in
V$, and $\sigma_1, \ldots, \sigma_{s-1} \in \{\pm 1\}$:
\begin{equation*}
\begin{split}
Y^{\sigma_1, \ldots, \sigma_{s-1}}_{i_1, \ldots, i_{s-1}} e_{i_s} & = Y^{\sigma_1, \ldots, \sigma_{s-1},+1}_{i_1, \ldots, i_{s}} e_\bO,\\
Y^{\sigma_1, \ldots, \sigma_{s-1}}_{i_1, \ldots, i_{s-1}} (e_\bO -
e_{i_s})& = Y^{\sigma_1, \ldots, \sigma_{s-1},-1}_{i_1, \ldots, i_{s}}
e_\bO,
\end{split}
\end{equation*}
where $Y^{\sigma_1, \ldots, \sigma_{s-1}}_{i_1, \ldots, i_{s-1}} = Y$ for $s = 1$.
\item[(c)] For all $i_1, \ldots, i_t \in V$ and $\sigma_1, \ldots,
\sigma_{t-1} \in \{\pm 1\}$:
\begin{equation*}
\begin{split}
Y^{\sigma_1, \ldots, \sigma_{t-1}}_{i_1, \ldots, i_{t-1}} e_{i_t} \in K,\\ 
Y^{\sigma_1, \ldots, \sigma_{t-1}}_{i_1, \ldots, i_{t-1}} (e_\bO - e_{i_t}) \in K,
\end{split}
\end{equation*}
where $Y^{\sigma_1, \ldots, \sigma_{t-1}}_{i_1, \ldots, i_{t-1}} = Y$ for $t = 1$.
\end{enumerate}
\end{proposition}

The above formulation allows to estimate the cost of optimizing over
$N_+^t(K)$ in terms of $n = |V|$ and $t$. Set $h(n,t):\sum_{s=0}^{t-1} (2n)^s=2^{t-1}n^{t-1} + O(n^{t-2})$. The formulation
involves ${n\choose 2}h(n,t)$ variables, $h(n,t)$ matrices of order
$n+1$, and $(2n)^t$ conditions of type ``$x\in K$''.  Furthermore, it
turns out that for $1\le r\le s$, the $i_r$-th column of the matrix
$Y^{\sigma_1, \ldots, \sigma_s}_{i_1, \ldots, i_s}$ is identically
zero whenever $\sigma_r=-1$ and it is equal to the $\bO$-th column
whenever $\sigma_r=1$. Thus, $Y^{\sigma_1, \ldots,
\sigma_s}_{i_1,\ldots,i_s}$ can be assumed to have order $n-s+1$.

\subsection{The Lasserre hierarchy}
\label{secL}

In this section we recall some basic facts about the Lasserre
construction, applied to the case when $K$ is of the form
\eqref{coneK}; for more information we refer to \cite{Las01b} and
\cite{Lau03}. The Lasserre hierarchy involves moment matrices: A
matrix whose rows and columns are indexed by a subset $\AAA$ of
$\PP(V)$ is said to be a {\em moment matrix} if the $(I,J)$-th entry
depends only on the union $I\cup J$ (for all $I,J\in \AAA$).  In this
definition one may allow $\AAA$ to be a multiset, which corresponds to
repeated rows and columns in the moment matrix. For a non-negative
integer $t$ and a vector $y\in \oR^{\PP_{2t}(V)}$, define the {\em
moment matrix of $y$ of order $t$} by $M_t(y):=(y_{I\cup
J})_{I,J\in\PP_t(V)}$. For a set $T$ and a vector $y\in \oR^{\PP(T)}$,
we write $M_T(y):=(y_{I\cup J})_{I,J\in\PP(T)}$; thus $M_n(y)=M_V(y)$
if $|V|=n$.

The following fact, observed in \cite{LS91, Lau03}, explains the
relevance of moment matrices to 0/1 polyhedra: For $x \in \oR^V$
define $\zeta_x:=(\prod_{i\in I}x_i)_{I\in \PP(V)}$. Then we have for
$y \in \oR^{\PP(V)}$
\begin{equation}\label{eqMVy}
\begin{split}
M_V(y)\succeq 0 & \Longleftrightarrow y\in \oR_+\left\{\zeta_x\mid x\in \{0,1\}^V\right\}\\
& \Longleftrightarrow  \forall S\subseteq V: \sum_{S': S \subseteq S' \subseteq V} (-1)^{|S'\setminus S|} y_{S'} \ge 0.
\end{split} 
\end{equation}
In Lemma~\ref{lemsymmetry} we give an extension of this result.

Next we explain how to encode the linear constraints $a_l^Tx \ge 0$
describing $K$. Given $y\in \oR^{\PP_{2t}(V)}$ and $a\in
\oR^{\PP_1(V)}$, define the vector $ay\in\oR^{\PP_{2t-1}(V)}$ by
$(ay)_I:= a_\bO y_I+\sum_{i\in V}a_iy_{I\cup \{i\}}$ for
$I\in\PP_{2t-1}(V)$.  For $t \geq 1$ we define the {\em $t$-th iterate of
the Lasserre hierarchy} by
\begin{equation}\label{setQtK}
\begin{split}
  Q^t(K):=\{x\in\oR^{\PP_1(V)}\mid \; & \exists y \in \oR^{\PP_{2t}(V)} :\;  
y_\bO=x_\bO,\ y_i=x_i\ (i \in V),\\
&  M_t(y)\succeq 0,\
 \ M_{t-1}(a_ly)\succeq 0\ (l=1,\ldots,m)\}.
\end{split}
\end{equation}
It lies between $K$ and $C_K$ and $Q^{t+1}(K)\subseteq Q^t(K)$. The Lasserre hierarchy refines the
Lov\'asz-Schrijver hierarchy, since we have $Q^{t+1}(K)\subseteq
N_+(Q^{t}(K))$ which implies $Q^{t+1}(K) \subseteq N_+^t(K)$ and
$Q^{n+1}(K)=C_K$. The formulation \eqref{setQtK} involves $\sum_{i= 0}^{2t} \binom{n}{i} = O(n^{2t})$ variables, one matrix of order
$\sum_{i = 0}^t \binom{n}{i} = O(n^t)$ and $m$ matrices of order $\sum_{i=0}^{t-1} \binom{n}{i} = O(n^{t-1})$.

\subsection{A new block-diagonal hierarchy}
\label{secnew}

One drawback of the Lasserre hierarchy is that the computational cost
for optimizing over $Q^t(K)$ is considerably higher than the cost for
optimizing over $N_+^{t-1}(K)$. To define a more economical variation
of it, which still refines the Lov\'asz-Schrijver hierarchy, we
consider a suitable principal submatrix of the full matrix $M_{t}(y)$.

For a positive integer $t$ and a subset $T\subseteq V$ of
cardinality $t-1$, let $M(T;y)$ denote the principal submatrix of
$M_{t}(y)$ whose rows and columns are indexed by
\begin{equation}
\label{indexset}
\AAA(T):=\displaystyle\bigcup_{S\subseteq T}\AAA_S, \ \text{ where }
 \AAA_S:=\{S \} \cup \{S\cup\{i\} \mid i\in V\}.
\end{equation}
It will be convenient to consider $\AAA(T)$ as a multiset: We keep
possible repeated occurrences, e.g.\ $S$ and $S\cup\{i\}$ if $i\in
S$. So strictly speaking the matrix $M(T;y)$ is a principal submatrix
of $M_{t}(y)$ only after removing repeated rows and columns. We
consider multisets here because it simplifies the notation in
Lemma~\ref{lemsymmetry}.
Note that $M(\emptyset;y) = M_1(y)$, and observe that we only need to
know the components of $y$ indexed by $\PP_{t+1}(V)$, instead of
$\PP_{2t}(V)$ as in the Lasserre hierarchy, in order to define the matrices
$M(T;y)$ for all $T\in\PP_{=(t-1)}(V)$.

Define the first iterate of the {\em block-diagonal hierarchy} by
$L^1(K):=Q^1(K)$ and, for $t\ge 2$, define its $t$-th iterate by
\begin{equation*}
\begin{split}
L^{t}(K) := \{x\in\oR^{\PP_1(V)}\mid & \exists y \in \oR^{\PP_{t+1}(V)}:
 y_\bO=x_\bO,\; y_i=x_i\ (i\in V),\\
& M(T;y)\succeq 0\ (T\in \PP_{=(t-1)}(V)),\\
& M(T;a_ly)\succeq 0\ (T\in \PP_{=(t-2)}(V), l=1,\ldots,m)\}.
\end{split}
\end{equation*}
Since we used principal submatrices of the Lasserre hierarchy, we
obviously have that the Lasserre hierarchy refines the block-diagonal
hierarchy. As we see in Section~\ref{ssec:comparisons} the
block-diagonal hierarchy still refines the Lov\'asz-Schrijver
hierarchy.

Next we give a more compact formulation for the set $L^t(K)$, based on
the fact that the matrix $M(T;y)$ has a special block structure which
can be exploited to block-diagonalize it. This property justifies the
name ``block-diagonal hierarchy''.

For a subset $S$ of $T$, let $A_S(y)$ denote the principal submatrix
of $M(T;y)$ indexed by the set $\AAA_S$, which is defined in \eqref{indexset}. It is a
$(n+1) \times (n+1)$ matrix lying in $\MMV$ with entries
\begin{equation*}
A_S(y)_{\bO,\bO}=y_S, \quad A_S(y)_{\bO,i}=y_{S\cup\{i\}}, \quad
A_S(y)_{i,j}=y_{S\cup\{i,j\}} \quad (i,j \in V).
\end{equation*}
The submatrix $M(T;y)[S,S']$ of $M(T;y)$ with row indices in $\AAA_S$
and column indices in $\AAA_{S'}$ depends only on $S\cup S'$:
$M(T;y)[S,S']=A_{S\cup S'}(y)$.

\begin{lemma}
\label{lemsymmetry}
The matrix $M(T;y)$ is positive semidefinite if and only if for all
subsets $S$ of $T$ the matrix
\begin{equation}\label{matAST}
A(S,T)(y) := \sum\limits_{S': S\subseteq
S'\subseteq T} (-1)^{|S'\setminus S|} A_{S'}(y)
\end{equation}
is positive semidefinite.
\end{lemma}

\begin{proof}
The proof is a ``block-matrix version'' of the one of \eqref{eqMVy} in
\cite{Lau03}. Define the block-matrix $Z$ indexed by $\AAA(T)$, whose
$(S,S')$-th block is the identity matrix $\bI$ of order $n+1$ if
$S\subseteq S'$ and the zero matrix otherwise. Its inverse is the
block matrix whose $(S,S')$-th block is $(-1)^{|S'\setminus S|}\bI$ if
$S\subseteq S'$ and the zero matrix otherwise.  Define the block
diagonal matrix $D$ with diagonal blocks $A(S,T)(y)$ for $S\subseteq
T$.  Direct verification shows that $M(T;y)=ZDZ^T$.  Therefore,
\begin{equation*}
M(T,y)\succeq 0\Longleftrightarrow D\succeq 0 \Longleftrightarrow
\forall S \subseteq T: A(S,T)(y)\succeq 0.\qedhere
\end{equation*}
\end{proof}

\begin{example} For $T=\{1,2\}$, 
$\AAA(T)=\AAA_\bO \cup \AAA_1\cup\AAA_2 \cup \AAA_{12}$ and
  \begin{align*}
M(T;y)=\left(\begin{matrix}
  A_\bO & A_{1} & A_{2} & A_{12}\cr
 A_1 & A_1 & A_{12} & A_{12}\cr
 A_2 & A_{12} & A_2 & A_{12}\cr
 A_{12} & A_{12} & A_{12} & A_{12}
\end{matrix}\right)
\succeq 0 \Longleftrightarrow 
\left\{\begin{array}{r}
A_\bO-A_1-A_2+A_{12}\succeq 0\\
A_1-A_{12}\succeq 0\\
A_2-A_{12}\succeq 0\\
A_{12}\succeq 0
\end{array}\right.
\end{align*}
where we wrote $A_S$ instead of $A_S(y)$.
\end{example}

Hence, in the formulation of $L^t(K)$, each condition $M(T;y)\succeq
0$, which involves one matrix of order $2^{t-1}(n+1)$, can be replaced
by the $2^{t-1}$ conditions $A(S,T)(y)\succeq 0$, each involving a
matrix of order $n+1$.  Similarly, the condition $M(T;a_ly)\succeq 0$
can be replaced by the $2^{t-2}$ conditions $A(S,T)(a_ly)\succeq 0$,
each involving a matrix of order $n+1$.

\subsection{A variation of the block-diagonal hierarchy}
\label{secvariation}

The next lemma deals with other possible ways of encoding the linear
conditions defining the set $K$. It motivates our second variation
$\tilde L^t(K)$. It turns out that it has an explicit link to the
Lov\'asz-Schrijver hierarchy. A main advantage of $\tL^t(K)$ over
$L^t(K)$ is that we do not need an explicit linear description of the
set $K$ in order to be able to define $\tL^t(K)$.  Hence $\tL^t(K)$
enjoys the same complexity property as $N_+^t(K)$: If one can optimize
in polynomial time over $K$ then the same holds for $\tL^t(K)$ for any
fixed $t$.

\begin{lemma}
\label{lem:linear conditions}
Let $t\ge 1$, $y\in \oR^{\PP_{t+1}(V)}$, $K$ be as in \eqref{coneK}
and $A(S,T)(y)$ be as in (\ref{matAST}).
Then, the following two assertions are equivalent:
\begin{enumerate}
\item[(a)] For all $T \in \PP_{=(t-1)}(V)$, $S \subseteq T$, $i \in
V$:
\begin{equation*}
A(S,T)(y)e_i \in K, \quad A(S,T)(y)(e_\bO-e_i) \in K.
\end{equation*}
\item[(b)] For all $T\in \PP_{\scriptscriptstyle=t}(V)$, $l=1,\ldots,m$:
\begin{equation*}
M_T(a_ly)\succeq 0.
\end{equation*}
\end{enumerate}
\end{lemma}

\begin{proof}
Using the identities
\begin{equation*}
a_l^TA_S(y)e_0=(a_ly)_S, \quad
a_l^TA_S(y)e_i=(a_ly)_{S\cup\{i\}},
\end{equation*}
the conditions $A(S,T)(y)e_i \in K,$ $A(S,T)(y)(e_\bO-e_i) \in K$ can
be rewritten as
\begin{equation*}
\begin{split}
&\sum_{S' : S\subseteq S'\subseteq T}(-1)^{|S'\setminus S|} (a_ly)_{S'\cup\{i\}}\ge 0 \ \
\ (l=1,\ldots,m),\\
&\sum_{S' : S\subseteq S'\subseteq T}(-1)^{|S'\setminus S|} 
((a_ly)_{S'}-(a_ly)_{S'\cup\{i\}})\ge 0 \ \ \ (l=1,\ldots,m).
\end{split}
\end{equation*}
On the other hand, using \eqref{eqMVy}, $M_T(a_ly)\succeq 0$ is equivalent to
\begin{equation*}
\sum_{S' : S\subseteq S'\subseteq T} (-1)^{|S'\setminus S|} (a_ly)_{S'} \ge 0 \ \ \ (S\subseteq T).
\end{equation*}
From this one can verify the equivalence of (a) and (b).
\end{proof}

Observe that for $t=1$ property (a) is equivalent to
$A_\bO(y)e_i,A_\bO(y)(e_0-e_i)\in K$ for all $i \in V$.  Combined with
the condition $A_\bO(y)\succeq 0$, this characterizes membership in
the set $N_+(K)$.  

This motivates replacing in the definition of $L^t(K)$ the condition
``$M(T;a_ly)\succeq 0$ for all $T\in \PP_{=(t-2)}(V)$'' by property
(a): For $t\ge 1$ define
\begin{equation*}
\begin{split}
\tL^t(K)
:=\{x\in\oR^{\PP_1(V)}\mid & \exists y \in \oR^{\PP_{t+1}(V)}:
 y_\bO=x_\bO, y_i=x_i\  (i\in V),\\
& M(T;y)\succeq 0 \ ( T\in \PP_{=(t-1)}(V)),\\
 & A(S,T)(y)e_i \in K,\; A(S,T)(y)(e_\bO-e_i) \in K\\
& \qquad (T \in \PP_{=(t-1)}(V), S \subseteq T, i \in V)\}.
\end{split}
\end{equation*}

\subsection{Comparisons}
\label{ssec:comparisons}

Another advantage is that $\tL^t(K)$ can be directly compared to the
Lov\'asz-Schrijver hierarchy $N_+^t(K)$. The next proposition shows
that our second variation refines the Lov\'asz-Schrijver hierarchy.

\begin{proposition}\label{lemtildeL}
We have $\tL^1(K)=N_+(K)$ and $\tL^t(K)\subseteq N_+^t(K)$ for $t\ge
2$.
\end{proposition}

\begin{proof}
As noted above we have $\tL^1(K)=N_+(K)$. Now let $t\ge 2$ and $x\in
\tL^t(K)$.  Thus, there is a $y\in \oR^{\PP_{t+1}(V)}$ which satisfies 
$ y_\bO=x_\bO,$ $y_i=x_i$ ($i\in V$), and 
$M(T;y)\succeq 0$ or, equivalently, $A(S,T)(y)\succeq 0$ for all 
$S\subseteq T\subseteq V$ with $|T|=t-1$.
Moreover property (a) of Lemma~\ref{lem:linear conditions} holds. 
Set $Y := M_1(y)$.  Then $x=Y e_\bO$ and $Y\in
\MMV$.  Given $1\le s\le t-1$, and $i_1,\ldots,i_s\in V$, and $\sigma\in \{\pm
1\}^s$, consider the multisets $T=\{i_1,\ldots,i_s\}$, $S=\{i_r \mid
r=1,\ldots,s, \sigma_r=1\}\subseteq T$, and 
define $Y^{\sigma_1,\ldots,\sigma_s}_{i_1, \ldots, i_s}:= A(S,T)(y)$. 
Here we extend the definition
of $A(S,T)(y)$ in (\ref{matAST})
to the case when $S$ and $T$ are multisets by taking the summation over all multisets $S'$ lying between $S$ and $T$; moreover,
 when $S'$ is a multiset with $S''$ as underlying set,
 we let $A_{S'}(y):=A_{S''}(y)$.
Now one can verify that the conditions from Proposition \ref{propeqNt} hold, which implies $x\in N_+^t(K)$.
\end{proof}

As one can see from the above proof, the main difference between $\tL^t(K)$
and $N_+^t(K)$ is that the matrices
$Y^{\sigma_1,\ldots,\sigma_s}_{i_1, \ldots, i_s}$ share many common
entries in the definition of $\tL^t(K)$. As a consequence, one can
describe the set $\tL^t(K)$ with less variables compared to
$N_+^t(K)$.  In Table~1 we compare the complexity
of the formulations for $\tL^t(K)$ and $N_+^t(K)$.  In both cases one
has a semidefinite programming formulation involving a number of
matrices of size $n+1$ required to be positive semidefinite and a
number of conditions of the type ``$x\in K$''.

{\small
\begin{table}[htb]
\begin{center}
\begin{tabular}{c|c|c}
  & $\tL^t(K)$  & $N_+^t(K)$   \\
\hline
\# variables  &  $\sum_{i=0}^{t+1} {n\choose i} $ &  ${n\choose 2}\sum_{i=0}^t (2n)^i $ \\
& $=\frac{1}{(t+1)!}n^{t+1} +O(n^t)$   &  $=2^{t-2}n^{t+1}+O(n^t)$ \\
\hline
\# matrices & ${n\choose t-1}2^{t-1}$ & $\sum_{i=0}^t (2n)^i$ \\
of order $n+1$  & $=\frac{2^{t-1}}{(t-1)!} n^{t-1} +O(n^{t-2})$ & 
$=2^{t-1} n^{t-1}+O(n^{t-2})$ \\
\hline
\# conditions & $2^t{n\choose t-1}$  & $2^tn^t$  \\
``$x\in K$'' & $ =\frac{2^t}{(t-1)!}n^{t-1}+O(n^{t-2})$   & \\
\hline
\end{tabular}
\end{center}
\caption{Complexity comparison for $\tL^t(K)$ and $N_+^t(K)$.}
\end{table}
\vspace*{-0.7cm}
}

Also, as already stated in Section~\ref{secnew}, the block-diagonal
hierarchy refines the Lov\'asz-Schrijver hierarchy. This can be seen
by comparing $L^{t+1}(K)$ with the second variation $\tL^t(K)$.

\begin{proposition}
For $t \geq 1$ the inclusion $L^{t+1}(K)\subseteq \tL^t(K)$ holds.
\end{proposition}

\begin{proof}
This follows directly from the definitions, after noting that, for
$|T|=t$, the index set of $M_T(y)$ is contained in the index set of
$M(T\setminus\{i\};y)$, where $i$ is any element of $T$.
\end{proof}

\section{Application to the stable set problem}
\label{secstable}

In this section we apply the new hierarchies to the stable set
problem. Let $G = (V, E)$ be a graph. A subset $S \subseteq V$ is
called a stable set if none of its vertices are adjacent.  The
incidence vector of $S$ is $\chi^S \in \{0,1\}^V$ with $\chi^S(i) = 1$
iff $i \in S$. The stability number $\alpha(G)$ is the maximum
cardinality of a stable set. By $\SSS_G$ we denote the set of all
stable sets of $G$. Then the stable set polytope is
\begin{equation*}
\STAB(G) := \conv\{\chi^S \mid S \in\SSS_G\},
\end{equation*}
and the corresponding cone is
\begin{equation*}
\ST(G):=\oR_+\left\{ \left(\begin{smallmatrix} 1 \\ \chi^S
\end{smallmatrix}\right) \mid S \in \SSS_G\right\}.
\end{equation*}
A linear relaxation of $\ST(G)$ is the fractional stable set cone
\begin{equation*}
\FR(G):=\{x\in \oR^{\PP_1(V)} \mid x_i\ge 0\ ( i\in V),\; x_i + x_j \leq x_\bO\ ( \{i,j\} \in E)\}.
\end{equation*}
A semidefinite relaxation of $\ST(G)$ is the theta body
\begin{equation*}
\TTH(G) := \{Ye_\bO \in \oR^{\PP_1(V)} \mid Y\in \MMV,
Y_{ij}=0\ ( \{i,j\} \in E)\},
\end{equation*}
which is contained in $\FR(G)$. Maximizing the linear function
$\sum_{i\in V}x_i$ over the theta body $\TTH(G)$ intersected with the
hyperplane $x_\bO = 1$ equals the Lov\'asz theta function $\vartheta(G)$
introduced by Lov\'asz in \cite{Lo79}. For details about these
relaxations and the stable set problem we refer e.g.\ to \cite{LR05}
and \cite{Sch03}.

In \cite[Lemma 20]{Lau03} it was shown that when constructing the
Lasserre hierarchy for $\FR(G)$ one can considerably simplify the
formulation. One can replace the condition ``$M_{t-1}(a_l y) \succeq
0$'', where $a_l$ runs through all linear inequalities defining
$\FR(G)$, by the simpler equalities $y_{ij} = 0$, where $\{i,j\} \in
E$, the so-called edge equalities. We want to apply the same
simplification to the definition of $L^t(\FR(G))$ and define another
variant $L^t(G)$ of it. However, in contrast to the Lasserre
hierarchy, this simplification weakens the block-diagonal hierarchy a
little bit since we can only claim the inclusion $L^t(\FR(G))\subseteq
L^t(G)$. Nevertheless the new variant $L^t(G)$ still refines the
Lov\'asz-Schrijver hierarchy, as $L^t(G)\subseteq N^{t-1}_+(\TTH(G))$
follows from Proposition~\ref{lemLtG} below combined with Proposition
\ref{lemtildeL}. We define
\begin{equation*}
\begin{split}
L^t(G) := \{x\in\oR^{\PP_1(V)}\mid & \exists y \in \oR^{\PP_{t+1}(V)} :
 y_\bO=x_\bO,\; y_i=x_i\ (i\in V),\\
& M(T;y)\succeq 0\ (T\in \PP_{=(t-1)}(V)),\\
& y_{ij}=0\  (\{i,j\} \in E)\}.
\end{split}
\end{equation*}
Thus, $L^1(G)=\TTH(G)$ and one can easily verify the inclusions
$\tL^t(\FR(G)) \subseteq L^t(G)$ when $t \geq 1$ and
$L^t(\FR(G))\subseteq L^t(G)$ when $t \geq 2$.  Maximizing the
objective function $\sum_{i\in V}x_i$ over $L^2(G)$ intersected with
the hyperplane $x_\bO = 1$ coincides with the parameter $\ell(G)$
considered in \cite{KP07,GL1,GL2,Lau07}.

The next lemma says that the edge conditions in the definition of
$L^t(G)$ imply that all variables indexed by non-stable sets are
identically $0$.

\begin{lemma}
\label{lemedge}
Let $y \in \oR^{\PP_{t+1}(V)}$ satisfy the
conditions in the definition of $L^t(G)$. Then $y_I=0$ for any subset $I\subseteq V$ with $|I|\le t+1$ and containing an edge.
\end{lemma}

\begin{proof}
For $|I|=2$ the statement is nothing else but the edge
equalities. Assume that $|I| \geq 3$, let $i, j\in I$ be adjacent
vertices, and let $k$ be another vertex in $I$. Define $T := I
\setminus \{i,j\}$. The matrix $M(T;y)$ is positive semidefinite and
the sets $\{i,j\}$ and $T \cup \{k\}$ occur in the index set
$\AAA(T)$. As the $(ij,ij)$-th entry of $M(T;y)$ is $y_{ij} = 0$, we
have by the positive semidefiniteness of $M(T;y)$ that its $(ij, T
\cup \{k\})$-th entry is $0$ as well and the statement of the lemma
follows.
\end{proof}

\begin{proposition} \label{lemLtG}
We have the inclusion $L^{t+1}(G)
\subseteq \tL^t(\TTH(G))$ for $t \geq 1$.
\end{proposition}

\begin{proof}
Assume that $y\in \oR^{\PP_{t+2}(V)}$ satisfies the conditions of the
definition of $L^{t+1}(G)$. In the following we show that the vector
consisting of the first $n+1$ coordinates of $y$ belongs to
$\tL^t(\TTH(G))$.

Fix $T\in \PP_{t-1}(V)$, $S\subseteq T$ and $k\in V$. We show that
$A(S,T)(y)e_k \in \TTH(G)$ and $A(S,T)(y)(e_\bO - e_k) \in
\TTH(G)$. For this we construct matrices $Y^k$ and $Z^k$ in $\MMV$
satisfying $Y^k_{ij} = Z^k_{ij} = 0$ when $i$ and $j$ are adjacent,
and satisfying
\begin{equation*}
Y^k e_\bO = A(S,T)(y)e_k,\quad Z^ke_\bO = A(S,T)(y)(e_\bO - e_k).
\end{equation*}
We distinguish between three cases.
\begin{enumerate}
\item $k \in S$:
Then $A(S,T)(y)e_\bO = A(S,T)(y)e_k$ and define $Y^k := A(S,T)(y)$, $Z^k := 0$.
\item $k \in T \setminus S$:
Then $A(S,T)(y)e_k = 0$ and define $Y^k := 0$, $Z^k := A(S,T)(y)$.
\item $k \in V \setminus T$:
Then we define $Y^k = A(S \cup \{k\}, T \cup \{k\})(y)$ and $Z^k:A(S,T\cup\{k\})(y)= A(S,T)(y)-A(S\cup\{k\},T\cup\{k\})(y)$.
\end{enumerate}
In all cases we see by Lemma \ref{lemsymmetry}, \ref{lemedge} that $Y^k,$  $Z^k$ satisfy the desired conditions.
\qedhere
\end{proof}

We summarize the inclusion relations between the various relaxations:
\begin{equation*}
\ST(G) \subseteq Q^t(\FR(G)) \subseteq L^t(G) \subseteq \tL^{t-1}(\TTH(G)) \subseteq N_+^{t-1}(\TTH(G)).
\end{equation*}
Moreover, $N^{t-1}_+(\TTH(G)) = \ST(G)$ holds for $t \geq \alpha(G)$
(see \cite{G} for a proof).

\section{Experimental results} \label{secPaley}

In this section we present some computational results for Paley graphs.

Let $\mathbb F_q$ be the finite field with prime power $q$ which is
congruent to $1$ modulo $4$; then $-1$ is a square in $\mathbb
F_q$. The {\it Paley graph} $P_q$ has $\mathbb F_q$ as vertex set and
two distinct elements $u,v\in \mathbb F_q$ are adjacent if $u-v$ is a square
in $\mathbb F_q$. The Paley graph is isomorphic to its complementary
graph, it is a strongly regular graph and its automorphism group acts
doubly-transitive on the vertices.  It is known (\cite[Theorem
8]{Lo79}) that $\vartheta(G)\vartheta(\ol G)=|V(G)|$ when $G$ is a
vertex-transitive graph, where $\ol G$ denotes the complementary graph
of $G$. Since the Paley graph $P_q$ is vertex-transitive and
isomorphic to its complementary graph, we have $\vartheta
(P_q)=\vartheta(\ol P_q)= \sqrt q$ (cf.\ \cite[Theorem
13.14]{Bol}). J.B.~Shearer (\cite{Sh96}) has computed $\alpha(P_q)$
for all primes $q\le 7000$.  For more information about $P_q$ we refer
e.g.\ to \cite[Chapter 13.2]{Bol}.

In order to illustrate the quality of the new relaxations $L^t(P_q)$,
we have computed the bounds obtained by maximizing $\sum_{v\in
  V(P_q)}x_v$ over the sets $L^t(P_q)$ (for $t=2,3$) and $N_+(\TTH(P_q))$ 
 intersected with $x_\bO = 1$. 
The results are given in Table~2. There we consider all
primes $q$ congruent to $1$ modulo $4$ between $61$ and $337$, as well
as a few larger values of $q$ up to $809$.  We have chosen the Payley
graph here because its automorphism group acts doubly-transitive on
the vertex set and so our formulation for $L^t(P_q)$ ($t\le 3$) and $N_+(\TTH(P_q))$ 
considerably simplifies.  (See  \cite[Chapter 6.1]{G} for  implementation details.)
For instance, optimization over
$L^3(P_{809})$ (resp., $L^2(P_{809})$, $N_+(\TTH(P_{809}))$)
can be formulated via an SDP with $876$ (resp., 36, 812) variables and
with four matrices with sizes 808, 808, 404 and 202 (resp.,
two matrices with sizes 809 and 405, three matrices with sizes 810, 810 and 809).
For the computations we used
the program CSDP \cite{Bor}.  Experiments were conducted on a single
machine with an 
Intel(R) Pentium(R) processor, 3Ghz and 1GB of RAM.
To compute the bounds from Table 2 we needed less
than a minute when $q \leq 100$ and,
for the largest instance $P_{809}$, around 45 minutes
for $L^3(P_{809})$, 31 minutes for $L^2(P_{809})$ and 4.5 hours for $N_+(\TTH(P_{809}))$.
Thus as expected the relaxation $L^2(G)$ gives a sharper bound than $N_+(\TTH(G))$, 
however at a much smaller computational cost.

Finally note that one can strengthen the relaxation $L^t(G)$ by adding the non-negativity constraints $y\ge 0$. However this only gives a  marginal improvement for
Paley graphs, as the bounds differ only in decimals.
{\small
\begin{table}[htb]
\begin{center}
\begin{tabular}{c|c|c|c|c|c}
$q$ 
& $\displaystyle{{L^1(P_q)=\TTH(P_q)}\atop {\vartheta(P_q)=\sqrt q}}$
& $N_+(\TTH (P_q))$
& $L^2(P_q)$ 
& $L^3(P_q)$ 
& $\alpha(P_q)$\\
\hline
61 & 7.810 & 5.901 & 5.465 & 5.035 & 5\\
73 & 8.544 & 6.377 & 5.973 & 5.132 & 5\\
89 & 9.434 & 7.155 & 6.304 & 5.391 & 5\\
97 & 9.849 & 7.948 & 7.398 & 6.596 & 6\\
101 & 10.050 & 7.290 & 6.611 & 5.496 & 5\\
109 & 10.440 & 8.007 & 7.366 & 6.578 & 6\\
113 & 10.630 & 8.330 & 7.599 & 7.009 & 7\\
137 & 11.705 & 8.829 & 8.200 & 7.047 & 7\\
149 & 12.207 & 9.188 & 8.231 & 7.136 & 7\\
157 & 12.530 & 9.695 & 8.707 & 7.485 & 7\\
173 & 13.153 & 10.316 & 9.426 & 8.062 & 8 \\
181 & 13.454 & 10.324 & 9.112 & 7.606 & 7 \\
193 & 13.892 & 10.506 & 9.210 & 7.651 & 7 \\
197 & 14.036 & 10.652 & 9.226 & 8.064 & 8 \\
229 & 15.133 & 11.659 & 10.290 & 9.076 & 9 \\
233 & 15.264 & 12.382 & 10.182 & 8.245 & 7 \\
241 & 15.524 & 11.595 & 9.891 & 8.275 & 7 \\
257 & 16.031 & 11.558 & 10.247 & 8.131 & 7 \\
269 & 16.401 & 12.307 & 10.624 & 8.778 & 8 \\
277 & 16.643 & 12.469 & 10.340 & 8.670 & 8 \\
281 & 16.763 & 11.902 & 10.605 & 8.397 & 7 \\
293 & 17.117 & 13.127 & 10.937 & 9.183 & 8 \\
313 & 17.692 & 13.128 & 11.630 & 9.458 & 8 \\
317 & 17.804 & 13.861 & 12.377 & 10.375 & 9 \\
337 & 18.358 & 13.724 & 11.658 & 9.464 & 9 \\
401 & 20.025 & 14.927 & 12.753 & 10.023 & 9\\
509 & 22.561 & 16.580 & 14.307 & 11.196 & 9\\
601 & 24.515 & 17.999 & 16.077 & 12.484 & 11\\
701 & 26.476 & 19.332 & 16.857 & 12.822 & 10\\
809 & 28.443 & 20.636 & 17.371 & 13.499 & 11\\
\hline
\end{tabular}
\label{tab:payley}
\medskip
\caption{Optimizing over $L^t(P_q)$ and $N_+(\TTH(P_q))$ for 
Paley graphs.}
\end{center}
\end{table}
}

\end{document}